\documentclass{amsart}
\usepackage{enumerate,amsmath,amssymb,amsfonts,bm}
\usepackage{frcursive,fullpage,graphicx,psfrag}
\usepackage{youngtab}
\usepackage{tikz}
\usetikzlibrary{shapes.geometric}
\usepackage{mathrsfs}\newcommand{\I}{\mathscr{P}}
 
\newcommand{\ZZ}{\mathbb{Z}}

\newcommand{\CC}{\mathbb{C}}

\newcommand{\aS}{\mathfrak{S}}
\newcommand{\ux}{\mathbf{x}}

\newtheorem{Theorem}{Theorem}

\theoremstyle{definition}
\newtheorem{Remark}{Remark}
\newcommand{\Weyl}[1]{\{{#1}\}}
\newcommand{\Specht}[1]{[{#1}]}
\newcommand{\tensor}{\smash{\textstyle\bigotimes}}
\newcommand{\Sym}{\mathsf{Sym}}

\begin{document}

\author{Abdelmalek Abdesselam}
\address{Abdelmalek Abdesselam,
Department of Mathematics,
P. O. Box 400137,
University of Virginia,
Charlottesville, VA 22904-4137, USA}
\email{malek@virginia.edu}

\author{Christian Ikenmeyer}
\address{Christian Ikenmeyer,
Department of Mathematics,
Mailstop 3368,
Texas A\&M University,
College Station, TX 77843-3368, USA}
\email{ciken@math.tamu.edu}

\author{Gordon Royle}
\address{Gordon Royle,
School of Mathematics and Statistics,
University of Western Australia,
35 Stirling Highway,
Crawley WA 6009, Australia}
\email{gordon.royle@uwa.edu.au}

\title{16,051 formulas for Ottaviani's invariant of cubic threefolds}

\begin{abstract}
We provide explicit combinatorial formulas for Ottaviani's degree 15 invariant
which detects cubics in 5 variables that are sums of 7 cubes.
Our approach is based on the chromatic properties of certain graphs and
relies on computer searches and calculations.  
\end{abstract}

\maketitle

\section{Introduction and main result}\label{introsec}
Throughout this article we will work over the field of complex numbers.
Let $v_d(\mathbb{P}^n) \subseteq \mathbb{P}^N$ denote the Veronese variety, where $N=\binom{n+d}{d}-1$.
The projective space $\mathbb{P}^N$ is the one built from the vector space of homogeneous polynomials $F$
of degree $d$ in $n+1$ variables. Inside this space, the Veronese variety corresponds to polynomials
which are of the form $F=L^d$ for some linear form $L$.
For a variety $X$ inside the projective space $\mathbb{P}^N$, let $\sigma_k(X)$ denote the $k$-th secant
variety of $X$, namely, the Zariski closure of the union of subspaces spanned by arbitrary collections of $k$ points in $X$.
In particular $\sigma_1(X)$ is $X$ itself.
Determining the codimension of such secant varieties is a difficult problem in general (see, e.g.,~\cite{Zak,ChiantiniC}). However, when
$X=v_d(\mathbb{P}^n)$ this question has been completely solved thanks to the celebrated theorem of Alexander and Hirschowitz.

\begin{Theorem}\cite{AlexanderH}
The codimension of $\sigma_{k}(v_d(\mathbb{P}^{n}))$
is given by the (naive) formula $\max\left\{N+1-(n+1)k,0\right\}$ except when:
\begin{enumerate}
\item
$d=2$ and $2\le k\le n$, for which the codimension is 
$\binom{n-k+2}{2}$;
\item
$d=4$, $n=2,3,4$ and $k=\frac{n(n+3)}{2}$ for which the codimension is $1$;
\item
$d=3$, $n=4$ and $k=7$ for which the codimension is $1$.
\end{enumerate} 
\end{Theorem}
 
A nice pedagogical introduction to this theorem can be found in~\cite[\S5.4 and Ch. 15]{Landsberg} (see also~\cite{BrambillaO}).
Of particular interest is the last and perhaps most involved exception in the previous theorem, namely $\sigma_7(v_3(\mathbb{P}^{4}))$.
This exception to the naive dimension-count formula was discovered around 1902 by Palatini~\cite{Palatini} and Richmond~\cite{Richmond}.
In 2009, Ottaviani~\cite{Ottaviani} proved that the hypersurface $\sigma_7(v_3(\mathbb{P}^{4}))$ was defined by an equation of degree 15
corresponding to the vanishing of an irreducible $SL_5$-invariant polynomial in the coefficients of the cubic $F$.
He also gave a formula for the cube of this invariant polynomial as the determinant of $45\times 45$ matrix with entries
that are linear in the coefficients of $F$. However, as far as we know, no formula for the invariant polynomial itself has
been given in the literature. The purpose of this article is to provide such an explicit formula.
In fact, we found 16,051 such formulas which are encoded by certain biregular bipartite graphs (these can be found in the accompanying
ancillary file \verb+GraphList.txt+). 
Only the one with highest symmetry is explicitly displayed in this article.
Note that the polynomial of interest is of degree 15 in 35 variables and could a priori have up to
$\binom{49}{15}=$1,575,580,702,584 monomials.
This number can be reduced, using invariance under the maximal torus of $SL_5$, so one need only consider isobaric
monomials. The number of such monomials is given (see, e.g.,~\cite[\S3.2]{OttavianiLect}) by the cofficient
of $x_1^9x_2^9x_3^9x_4^9 y^{15}$ in
\[
\prod_{\substack{i_1,\ldots,i_4\ge 0 \\ i_1+\cdots+i_4\le 3}}\frac{1}{1-x_{1}^{i_1}x_{2}^{i_2}x_{3}^{i_3}x_{4}^{i_4}y}\ .
\]
Equivalently, this number is that of equivalence classes under column permutation of $5\times 15$ arrays of nonnegative integers
with row sums equal to 9 and column sums equal to 3. Via this second characterization, we computed
this number using the constraint solver \textsc{Minion}~\cite{GentJM} available at \verb+http://minion.sourceforge.net/+
Forcing one representative in each class can be conveniently implemented in that software by adding constraints corresponding
to an imposed lexicographic ordering for the columns of the array. After about 551 seconds of computing time, the returned
number for isobaric monomials was
317,881,154 which still is an astronomically large number.
An explicit expansion of the sought polynomial is therefore out of the question. 
Instead we will give a compact combinatorial formula using the Feynman diagram calculus (FDC) explained in~\cite[\S2]{Abdesselam}.
This graphical formalism is equivalent to the symbolic method of classical invariant theory.
Our main finding is the following theorem.

\begin{Theorem}\label{mainthm}
The hypersurface $\sigma_7(v_3(\mathbb{P}^{4}))$ is defined by the equation $\I(F)=0$ where $\I(F)$ is the polynomial given
by
\[
\hspace{-2cm}\I(F)=\hspace{10cm}
\]
\[
\parbox{6cm}{\psfrag{a}{$\scriptstyle{14}$}
\psfrag{b}{$\scriptstyle{9}$}
\psfrag{c}{$\scriptstyle{10}$}
\psfrag{d}{$\scriptstyle{13}$}
\psfrag{e}{$\scriptstyle{1}$}
\psfrag{f}{$\scriptstyle{3}$}
\psfrag{g}{$\scriptstyle{2}$}
\psfrag{h}{$\scriptstyle{6}$}
\psfrag{i}{$\scriptstyle{11}$}
\psfrag{j}{$\scriptstyle{0}$}
\psfrag{k}{$\scriptstyle{12}$}
\psfrag{l}{$\scriptstyle{4}$}
\psfrag{m}{$\scriptstyle{5}$}
\psfrag{n}{$\scriptstyle{8}$}
\psfrag{o}{$\scriptstyle{7}$}
\hspace{-5cm}
\includegraphics[width=6cm]{FDCpicture.eps}}
\]

\vskip 5cm
\noindent
in the FDC notation explained in~\cite[\S2]{Abdesselam}.
\end{Theorem}

The graphical notation we used can be briefly described as follows.
Let our variable cubic $F(\ux)=F(x_1,\ldots,x_5)$ be written using multinomial coefficients
as
\[
F(\ux)=\sum_{|\alpha|=3}
\binom{|\alpha|}{\alpha}
\ f_\alpha \ux^\alpha
\]
where $\alpha=(\alpha_1,\ldots,\alpha_5)$ is a multiindex of nonnegative integers,
$|\alpha|=\sum_{i=1}^{5}\alpha_i$, $\alpha!=\prod_{i=1}^{5}\alpha_i!$, $\ux^\alpha=\prod_{i=1}^{5}x_i^{\alpha_i}$
and
\[
\binom{|\alpha|}{\alpha}
=\frac{|\alpha|!}{\alpha!}\ .
\]
Alternatively, using the identification of $F$ with a symmetric tensor, one can write
\[
F(\ux)=\sum_{i_1,i_2,i_3=1}^{5}
F_{i_1,i_2,i_3}\ x_{i_1} x_{i_2} x_{i_3}
\]
with $F_{i_1,i_2,i_3}=f_\alpha$ where the multiindex is defined by
\[
\alpha_i=|\{\nu\in\{1,2,3\}| i_\nu=i\}|\ .
\]
Namely, $\alpha_i$ counts how many times the value $i$ appears in the sequence of indices $(i_1,i_2,i_3)$.
The tensor entry $F_{i_1,i_2,i_3}$ is symmetric with respect to the permutation of the three indices.
We now introduce a graphical notation (an $F$-blob) for such a tensor entry:
\[
\parbox{1.2cm}{\psfrag{a}{$\scriptstyle{i_1}$}
\psfrag{b}{$\scriptstyle{i_2}$}
\psfrag{c}{$\scriptstyle{i_3}$}
\includegraphics[width=1.2cm]{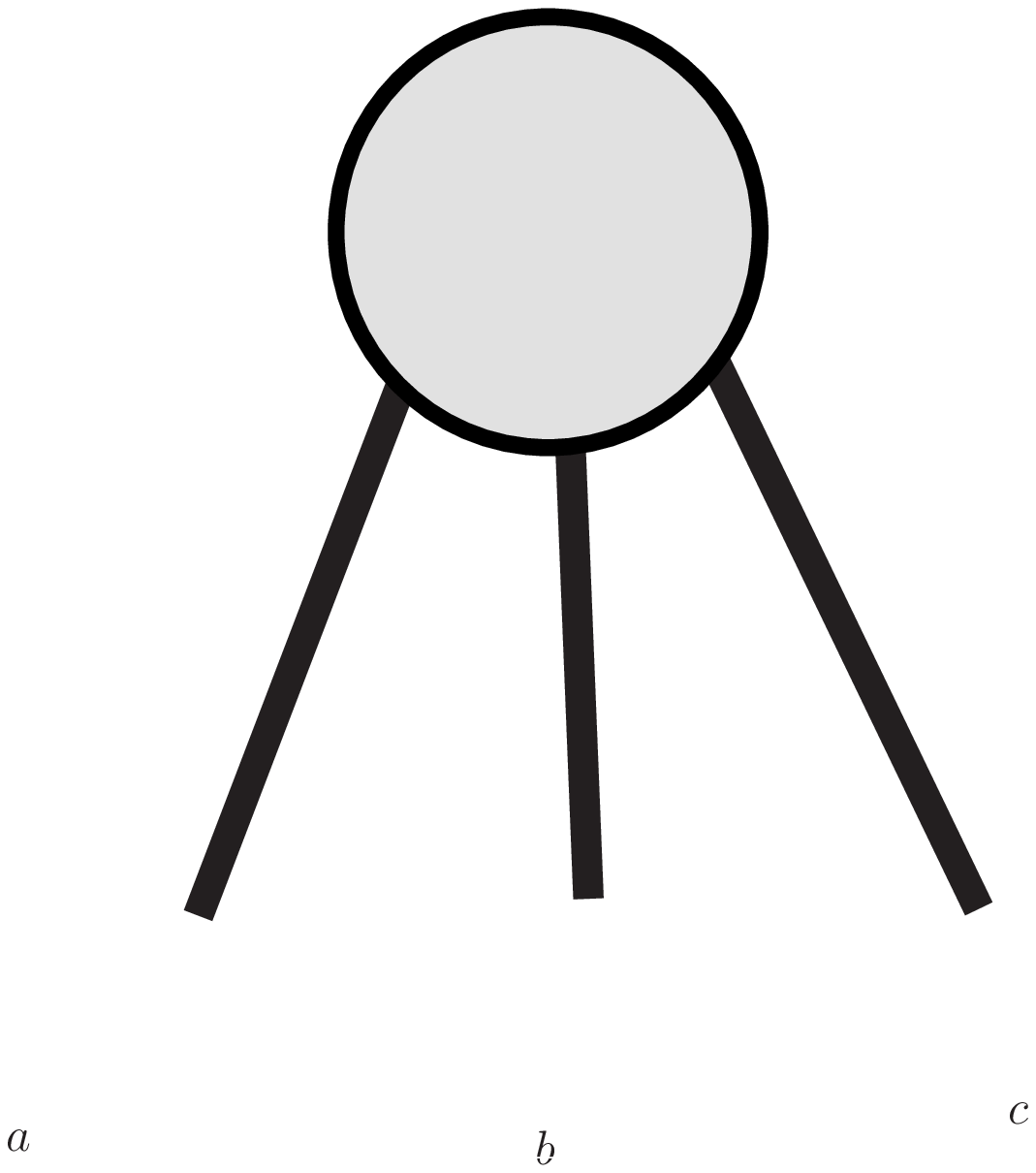}}
=F_{i_1,i_2,i_3}\ .
\]
Each of the three legs (or half-edges) of an $F$-blob therefore carries an index with values in $\{1,2,3,4,5\}$.
Likewise, we also define the completely antisymmetric tensor $\epsilon=(\epsilon_{i_1,\ldots,i_5})_{1\le i_1,\ldots,i_5\le 5}$ as follows.
If two indices are repeated we let $\epsilon_{i_1,\ldots,i_5}=0$. Otherwise we let $\epsilon_{i_1,\ldots,i_5}$ be equal to the
sign of the permutation $i_1,i_2,\ldots,i_5$ of $1,2,\ldots,5$.   
We similarly introduce a graphical notation (an $\epsilon$-fan) for this object
\[
\parbox{1.5cm}{\psfrag{a}{$\scriptstyle{i_1}$}
\psfrag{b}{$\scriptstyle{i_2}$}
\psfrag{c}{$\scriptstyle{i_3}$}
\psfrag{d}{$\scriptstyle{i_4}$}
\psfrag{e}{$\scriptstyle{i_5}$}
\includegraphics[width=1.5 cm]{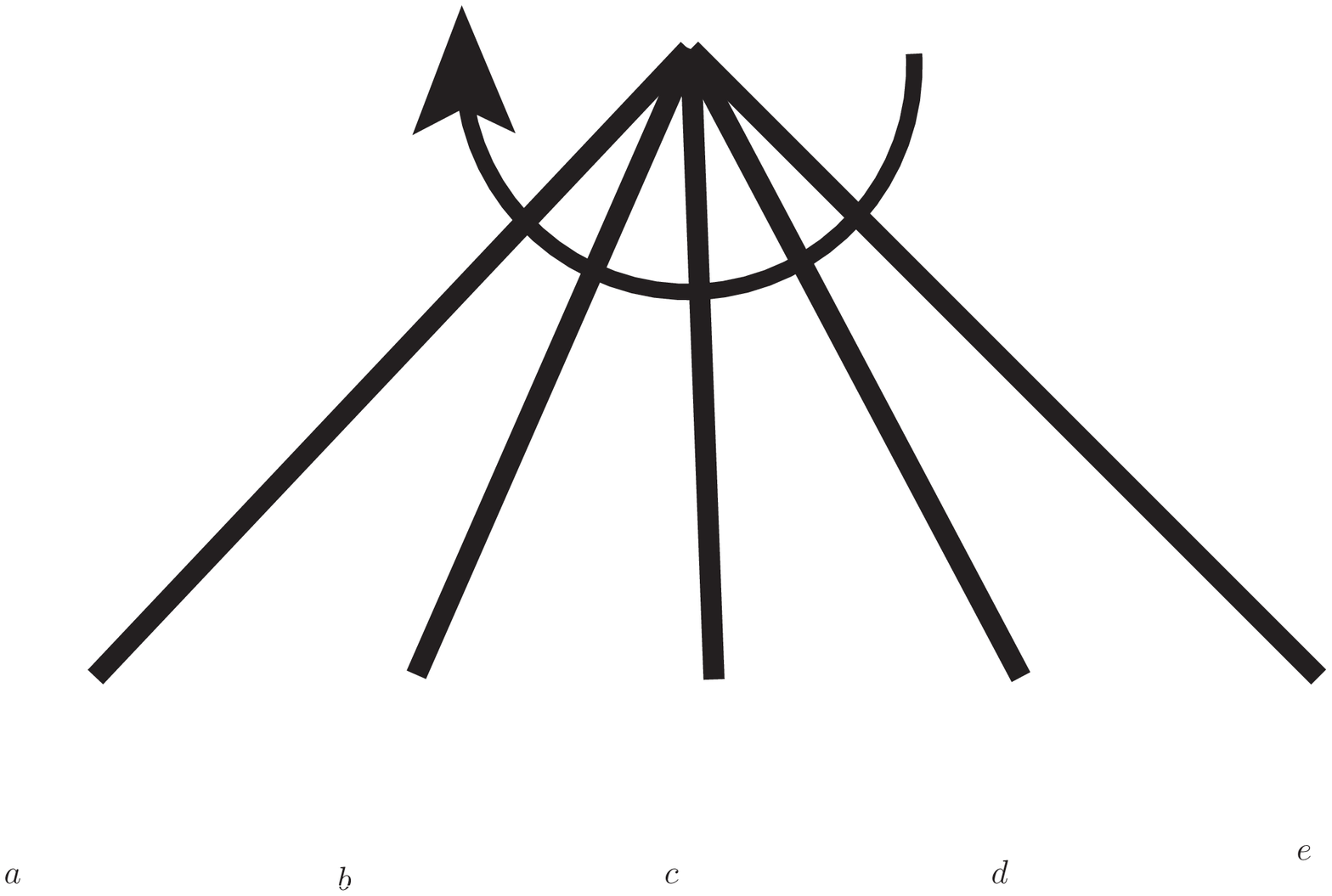}}
=\epsilon_{i_1,i_2,i_3,i_4,i_5}\ .
\]
Similarly, each of the five legs of an $\epsilon$-fan carries an index 
with values in $\{1,2,3,4,5\}$. Here the ordering of the legs is important in order to avoid sign ambiguities.
We will call the leg carrying the index $i_1$ the 1st leg, the one carrying $i_2$ the 2nd leg etc.
Thus, our curly arrow notation is the reverse of that introduced by Creutz in~\cite{Creutz1,Creutz2}.

To a diagram obtained by assembling $F$-blobs and $\epsilon$-fans,
such as the one in the statement of Theorem \ref{mainthm}, one associates a polynomial
in the coefficients $f_\alpha$ of the cubic $F$.
One assigns the same index to each pair of glued legs, then takes the product of corresponding tensor elements,
and finally sums over the indices in the set $\{1,2,3,4,5\}$.
This simply gives a compact notation for a complicated contraction of tensors (with summation over 45 indices).
The polynomial $\I(F)$ is homogeneous of degree 15 in the $f_\alpha$'s because it involves 15 $F$-blobs.

An alternate and equivalent description of such polynomials can also be given in terms of the symbolic method of classical invariant theory
as follows.
Let $a=(a_1,\ldots,a_5)$, $b=(b_1,\ldots,b_5)$, etc. denote vectors of independent variables.
Determinants obtained from such (row) vectors are denoted by:
\[
(abcde)=\left|
\begin{array}{ccccc}
a_1 & a_2 & a_3 & a_4 & a_5 \\
b_1 & b_2 & b_3 & b_4 & b_5 \\
c_1 & c_2 & c_3 & c_4 & c_5 \\
d_1 & d_2 & d_3 & d_4 & d_5 \\
e_1 & e_2 & e_3 & e_4 & e_5 
\end{array}
\right|
\]
and similarly for other combinations such as $(abfgl)$, etc.
Consider the product of determinants
\begin{equation}
\mathcal{S}=(bdcae)(gafih)(jknlm)(ojkbf)
(bdcgl)(mefih)(ojnch)(okndi)(glame)\ .
\label{symboliceq}
\end{equation}
The latter is a polynomial in $\CC[a_1,\ldots,o_5]$ (75 variables) which is separately homogeneous
of degree 3 in each letter $a, b,\ldots, o$.
Now consider the differential operator of order $45$
\[
\mathcal{D}=\frac{1}{3!}F\left(\frac{\partial}{\partial a_1},\ldots,\frac{\partial}{\partial a_5}\right)\ 
\frac{1}{3!}F\left(\frac{\partial}{\partial b_1},\ldots,\frac{\partial}{\partial b_5}\right)\ 
\cdots\ 
\frac{1}{3!}F\left(\frac{\partial}{\partial o_1},\ldots,\frac{\partial}{\partial o_5}\right)\ 
\]
where partial derivatives are substituted for the variables.
Then the result of applying $\mathcal{D}$ on the expression $\mathcal{S}$
is exactly equal to the polynomial $\I(F)$ from Theorem \ref{mainthm}. The labels from 0 to 14 provided
in the defining picture for $\I(F)$ correspond to the letters $a, b,\ldots, o$ numbered in the alphabet order from 0 to 14.

Consider disjoint finite sets $\mathcal{A}$ and $\mathcal{B}$ with respective cardinalities 15 and 9.
We will in fact take $\mathcal{A}=\{0,1,\ldots,14\}$.
Suppose one has a biregular bipartite graph on the vertices given by the partition formed by $\mathcal{A}$ and $\mathcal{B}$
where each vertex in $\mathcal{A}$ has degree 3 and each vertex in $\mathcal{B}$ has degree 5.
It is clear that the FDC or the symbolic method allow one to associate a polynomial as above to such a bipartite graph,
up to a sign ambiguity.
The latter can be removed by specifying an ordering of the edges at each vertex of $\mathcal{B}$.
Alternately, one can encode the polynomial by a hypergraph on the vertex set $\mathcal{A}$ made of 9 hyperedges (or 5-blocks)
which are subsets of 5 elements in $\mathcal{A}$. Such hypergraphs are also called 1-(15,5,3) designs or tactical configurations.
For instance, the previous hypergraph can be specified by the list
\[
1,3,2,0,4;\ 6,0,5,8,7;\ 9,10,13,11,12;\ 14,9,10,1,5;\ 1,3,2,6,11;\ 12,4,5,8,7;\ 14,9,13,2,7;
\]
\[
14,10,13,3,8;\ 6,11,0,12,4
\]
where 5-blocks have been separated by semicolons and extra spaces.
Via the alphabet numbering correspondence between the letters $a$, $b$,\ldots, $o$ and the numbers $0$, $1$,\ldots, $14$,
the above list exactly matches the symbolic expression in (\ref{symboliceq}).
The previous list was produced as the 21,238-th line of a computer-generated file.
In fact the computer produced the list
\[
0,1,2,3,4;\ 0,5,6,7,8;\ 9,10,11,12,13;\ 1,5,9,10,14;\ 1,2,3,6,11;\ 4,5,7,8,12;\ 2,7,9,13,14;
\]
\[
3,8,10,13,14;\ 0,4,6,11,12
\]
The reordering in each 5-block (for the first list)
is simply due to trying to make the picture\ in Theorem \ref{mainthm} prettier (avoiding unnecessary
twisting and crossing of $\epsilon$-fan legs).
The two lists give the same polynomial $\I(F)$ (with the same sign) since, out of the 9 reordering permutations, exactly 4 have negative sign.

Note that the ordering of $5$-blocks does not affect the resulting polynomial $\I$ (this amounts to moving the determinant factors
in the formula for $\mathcal{S}$).
The polynomial $\I$ is also unaffected by permuting the symbolic letters in the symbolic expression $\mathcal{S}$
because the differential operator $\mathcal{D}$ is symmetric in these
letters. As a result, in our search for suitable bipartite graphs we only need to exhaustively analyze their isomorphism classes
with respect to permutation of the vertices in $\mathcal{A}$ as well as permutation of vertices in $\mathcal{B}$.

One can also relate our constructions to modern representation theory and semistandard tableaux as follows.
By Schur-Weyl duality, the tensor space $\tensor^{45}\CC^5$ can be decomposed into
irreducibles as \mbox{$\bigoplus_{\lambda} \Specht\lambda \otimes \Weyl\lambda$},
where the sum is over all Young diagrams $\lambda$ with 45 boxes and at most 5 rows
and $\Specht\lambda$ denotes the irreducible $\aS_{45}$ representation of $\lambda$ (the Specht module)
and $\Weyl\lambda$ denotes the irreducible $GL_5$ representation of $\lambda$ (the Weyl module).
Invariants under $SL_5$ are irreducible representations of $GL_5$ that correspond to Young diagrams with five rows of length 9 each,
i.e., the Young diagram $(9^5)$.
The $SL_5$ invariant space is hence parametrized by the Specht module $\Specht{9^5}$.
A basis can be given by standard tableaux of shape $(9^5)$.
If we project this basis of the $SL_5$ invariant space of $\tensor^{45}\CC^5$ to the the $SL_5$ invariant space of $\tensor^{15}(\Sym^3\CC^5)$,
then we obtain the basis of $\tensor^{15}(\Sym^3\CC^5)$ which is given by semistandard tableaux of shape $(9^5)$ and content $(3^{15})$.
A tableau can be contructed from a graph as in Theorem~\ref{mainthm} by numbering the $F$-blobs and placing in each column the numbers
of the $F$-blobs that correspond to $\epsilon$-fans. Note that the order in each column and the order of the columns does not matter,
so that we can get different tableaux from the same graph.
One tableau corresponding to the graph in Theorem~\ref{mainthm} is
{
\def\A{10}
\def\B{11}
\def\C{12}
\def\D{13}
\def\E{14}
\def\F{15}
\begin{center}
\young(159114234,26\A525786,3\F\B9379\A\B,47\C\A68\D\D\C,\F8\D\E\B\C\E\E\F)
\end{center}
}
\noindent
(where we replaced 0 with 15).
This tableau is not semistandard, but can be expressed over the basis of semistandard tableaux using the so-called
straightening algorithm, which can be found for example in \cite[p. 110]{Fulton}.

The proof of Theorem \ref{mainthm} is an immediate consequence of the following two properties satisfied by our graph (or by any other
in our list of 16,051 graphs). Both of these properties are not obvious and have been checked by computer.

\smallskip
\noindent{\bf The chromatic property:} The graph does not have a proper coloring with 7 colors.

\smallskip
\noindent{\bf The nonvanishing property:} The associated polynomial $\I(F)$ is not identically zero.

\smallskip
By proper coloring with $q$ colors we mean an assignment of colors among $\{1,2,\ldots,q\}$ for every vertex of $\mathcal{A}$
such that for every 5-block the vertices in that 5-block are given 5 distinct colors. 
The chromatic property immediately implies that $\I(F)$ vanishes on the desired locus $\sigma_7(v_3(\mathbb{P}^4))$.
Indeed if $F$ is a sum of 7 cubes, $F=L_1^3+\cdots+L_7^3$, then the multilinear expansion of $\I(F)$ will
produce $7^{15}$ terms which are products of determinants as in (\ref{symboliceq}) where each letter has been replaced by a vector
of coefficients of some of the linear forms $L$, as indicated by the coloring.
For instance, the letter $g$ corresponds to the number or vertex 6 in $\mathcal{A}=\{0,1,\ldots,14\}$.
If that vertex has been assigned color $\nu\in\{1,\ldots,7\}$, then the vector of variables $g=(g_1,\ldots,g_5)$
is replaced by the vector of coefficients of the linear form $L_{\nu}$. A similar statement holds for the other letters, as is
most easily seen by expanding the differential operator $\mathcal{D}$ in terms of the linear forms.
Since no such 7-coloring is proper (by the chromatic property), there must be a determinant factor (5-block) which contains
repeated rows (repeated colors). As a result, all the terms of the previous multilinear expansion vanish and hence $\I(F)=0$.

\begin{Remark}\label{Sinvremark}
A similar phenomenon occurs for the Aronhold $S$-invariant of ternary cubics.
This invariant was introduced in~\cite{Aronhold1} where it was also shown to provide an equation for $\sigma_{3}(v_{3}(\mathbb{P}^2))$.
Similarly to (\ref{symboliceq}), Aronhold's invariant can be given an explicit symbolic expression as $(abc)(abd)(acd)(bcd)$
which was discovered in~\cite{Aronhold2}.
The vertices of the corresponding graph can be seen as those of a tetrahedron with the faces providing the hyperedges of $3$-blocks.
In this situation the chromatic property is trivial: if one colors the vertices with 3 colors only, then a face must have two
vertices of the same color. By the same reasoning as above, the corresponding degree 4 invariant vanishes on the locus of ternary
cubics which are sums of 3 cubes of linear forms (see also~\cite[p. 57]{OttavianiLect}).
\end{Remark}

In the graph search explained with more detail in \S\ref{graphsec}, we found 21,282 viable (isomorphism classes of) graphs
which satisfy the chromatic property. A polynomial $\I$ corresponding to such a graph, being of degree 15 by construction
and vanishing on $\sigma_7(v_3(\mathbb{P}^4))$ by the chromatic property, must be a constant multiple of the defining equation
of that variety. Indeed, Ottaviani~\cite{Ottaviani} showed that the latter is defined by an irreducible polynomial of degree 15.
Thus, all that remains to do is to check that this proportionality constant is nonzero, i.e., that the polynomial
given by the graph at hand is not identically zero.
This was checked by evaluation on a fixed random sum of 8 cubes.
Details of this procedure are given in \S\ref{nonvansec}.
Note that if our computation says the polynomial is nonzero, then we produced a computer proof of the latter statement.
The other direction holds as well: If our computation says the polynomial is zero, then we have a computer proof of this statement.
This follows from the fact that all polynomials under consideration are the same up to a constant factor and all are evaluated at the same point\footnote{Evaluating at the same point was suggested by C.~Krattenthaler. Moreover, in a previous version of this article, our results in the second direction
concerning the vanishing statement were only certified with a high probability but not absolute degree of certainty.
This is because, in order to speed up calculations, we used a mod-$p$ reduction. Recently we were able to improve
our computational result by working over $\ZZ$ and produce a computer proof of the identical vanishing property
for the 5,170 graphs for which this occurs.}.
Table \ref{tab:tallybyG} gives the count of graphs which pass this evaluation test or not, organized by the size of their symmetry group $G$.

\begin{table}
\begin{center}
\begin{tabular}{|c|c|c|c|}
\hline
$|G|$ & nonzero & zero & total \\ \hline 
32 & 0 & 1 & 1\\ \hline
16 & 0 & 2 & 2\\ \hline
12 & 1 & 0 & 1\\ \hline
8 & 8 & 6 & 14\\ \hline
6 & 1 & 0 & 1\\ \hline
4 & 58 & 51 & 109\\ \hline
3 & 0 & 1 & 1\\ \hline
1 or 2 & 15,983 & 5,170 & 21,153\\ \hline
\end{tabular}
\end{center}
\caption{Viable hypergraphs organized by the size of their automorphism group}
\label{tab:tallybyG}
\end{table}

To summarize the results of our search,
of the 21,282 graphs, 16,051 provide a nonzero polynomial $\I$ whereas 5,231 of them provide polynomials which are 
identically zero.
Here the symmetry group $G$ is the set of permutations of the vertex set $\mathcal{A}$ which conserve the hypergraph structure
(i.e., preserves the unordered set of 5-blocks which themselves are seen as unordered subsets of five elements in $\mathcal{A}$).
The graph chosen in the statement and picture for Theorem \ref{mainthm} is the unique graph with $|G|=12$.
Its symmetry group is generated by the permutations
\[
\tau=(2\ 3)(7\ 8)(9\ 10)
\]
and
\[
\omega=(1\ 7\ 3\ 5\ 2\ 8)(9\ 13\ 10)(4\ 6)(11\ 12)
\]
given as products of cycles. The corresponding relations are $\tau^2=\omega^6={\rm Id}$ and $\tau\omega\tau=\omega^{-1}$
which makes $G$ isomorphic to the dihedral group of order 12.

Finally, note that the relation between secants of Vernonese varieties and coloring properties of graphs which is
used in this article is different and rather more straightforward than the one featuring in~\cite{SturmfelsS}.
Also, a method similar to the one in this article has been used in~\cite{BI:13,HauensteinIL,Ikenmeyer}
regarding the border-rank problem.

\section{The graph search}\label{graphsec}

The bipartite incidence graph between the points $\mathcal A$ and hyperedges $\mathcal B$ is a biregular bipartite
graph with $15 + 9 = 24$ vertices. This is small enough that it is possible, indeed straightforward, to compute all
the possible graphs and then test which of them meet the chromatic and nonvanishing requirements described above.
Brendan McKay's program \textsc{genbg}, distributed with the \textsc{nauty Traces} package of
McKay and Piperno \cite{McKayPiperno} is a 
highly optimised implementation of an orderly algorithm for generation of classes of pairwise non-isomorphic bipartite graphs.
The command \verb+genbg -d3:5 -D3:5 15 9+
where the switches  {\tt -d} and {\tt -D} constrain the minimum and maximum degrees of the vertices in each part, constructs  2,553,049
candidate graphs in just under 30 minutes cpu time (3.4GHz Intel Core i7).

With no {\em a priori} restrictions, the hypergraphs may have {\em multiple hyperedges} (two or more identical hyperedges) or {\em repeated
vertices} (two vertices in the same set of hyperedges). There is no {\em a priori} reason to exclude the former situation, however
the latter leads to hypergraphs which are not viable for our purposes.
Indeed, if a hypergraph has repeated vertices the associated polynomial $\I$ must be identically zero.
This can be seen in the FDC as follows: exchange the $F$-blobs corresponding to the pair of repeated vertices simply by moving them
(which does not change $\I(F)$), then switch the $\epsilon$-fan legs where they are attached in order to recover the same diagram
one started with.
By antisymmetry of the $\epsilon$ tensors this produces a sign change of $(-1)^3$, so the polynomial is equal to its opposite and thus vanishes.
Of the 2,553,049 hypergraphs constructed there are 644,897 with no repeated vertices. The other ones are eliminated at this stage
of the search.

The chromatic property requires that the $15$-vertex ``collinearity graph'' with vertex set $\mathcal A$ and where $v$, $w$ are 
adjacent if they lie in a common hyperedge, have no proper $7$-coloring. (Equivalently we can view the edge set of the collinearity
graph as the union of nine $5$-cliques corresponding to the nine hyperedges).  Determining the chromatic number $\chi$ of a $15$-vertex graph is 
very easy by computer (a few seconds for all of the graphs).
The results, namely the number of graphs for given $\chi$, are shown in Table~\ref{tab:colouring}.

\begin{table}
\begin{center}
\begin{tabular}{|c|r|}
\hline
$\chi$ & Number\\
\hline
5 & 13,711\\
6 & 212,058\\
7 & 397,730\\
8 & 21,398\\
\hline
\end{tabular}
\end{center}
\caption{Chromatic numbers of the collinearity graphs}
\label{tab:colouring}
\end{table}

The easiest way for a graph to be $8$-chromatic is for it to possess an $8$-clique; however in this case,
the associated polynomial is identically
zero. More generally, a graph with $\chi=8$ will automatically give a zero polynomial if it is not  {\em vertex-critical}.
Here, vertex-critical means a graph such that $\chi(G) = 8$ and $\chi(G-v) < 8$ for all vertices $v \in V(G)$.
This property is also easy to see in the FDC language. If there is a vertex $v$ with $\chi(G-v)=8$
then remove the corresponding $F$-blob and assign fixed indices
$i$, $j$, $k$ to the $\epsilon$-fan legs left hanging (this is similar to the classical notion of evectant,
see e.g.~\cite[\S2.3]{DAndreaC}).
This produces a tensor $(T_{i,j,k}(F))_{1\le i,j,k\le 5}$ with entries that are degree 14 polynomials in the coefficients of the cubic $F$.
The same mutilinear expansion argument (for $F$ of the form $L_1^3+\cdots+L_7^3$) used earlier and the hypothesis $\chi(G-v)=8$
imply that each $T_{i,j,k}(F)$ vanishes on $\sigma_7(v_3(\mathbb{P}^4))$ which is of degree 15 by the result in~\cite{Ottaviani}.
Therefore $T_{i,j,k}(F)$ being of degree 14 must vanish identically in $F$ and for every choice of $i,j,k$.
Since the original invariant $\I$ is equal to $\sum_{i,j,k} T_{i,j,k}(F) F_{i,j,k}$, it will vanish too.  
Of the 21,398 remaining graphs, exactly 116 are not vertex-critical (112 which contain an $8$-clique and 2 which do not)
and need no further consideration.

In sum, we are left with 21,282 graphs which satisfy the chromatic property and
do not vanish for obvious reasons. These graphs are passed through into the nonvanishing check. 
We also note that none of them has repeated hyperedges.

\section{The nonvanishing check}\label{nonvansec}

We choose random vectors $L_1, \ldots, L_8 \in \ZZ^5$ and evaluate all polynomials at this random point.
As explained in \S\ref{introsec}, we need to compute $\sum_{c}\text{eval}(c)$,
where the sum is over all proper colorings $c$ of the graph with 8 colors and
$\text{eval}(c)$ is the corresponding product of determinants as described in \S\ref{introsec}.

Our fairly standard approach is as follows.
To list all proper colorings with 8 colors we use a depth-first backtracking algorithm,
in principle the same as is used when solving the famous eight queens puzzle with a recursive procedure.
Instead of checking all $8^{15}$ colorings, we color $F$-blobs one by one,
making sure that at each step no two $F$-blobs of the same color are connected
by an $\epsilon$-fan.
Whenever we reach a situation where there is an $F$-blob that cannot be colored without violating this condition,
we undo the last coloring step and choose
the next color for that vertex.
The order in which we color the $F$-blobs has a large influence on the running time of our algorithm,
so we use the greedy low-width heuristic, which is as follows:
at every stage, each uncolored $F$-blob has a set of colors that could be assigned to it such that no two $F$-blobs of the same
color are connected by an $\epsilon$-fan.
From all uncolored $F$-blobs we choose (in an arbitrary way) one that has a minimal number of possible colors. 
There are 38,707,200 proper colorings $c$ for the graph shown in Theorem \ref{mainthm}
and listing them while summing over their evaluations $\text{eval}(c)$
takes half an hour on a laptop computer.

To perform the similar calculation for all 21,282 graphs, we used the Texas A\&M University Brazos HPC cluster computer.

We took the evaluation results for all graphs, divided by their gcd, and listed this number in the ancillary file attached to this document.
The absolute values that appear as results are 1, $\ldots$, 30, 32, 33, 35, 36, 39, 40, 42, 45, 48, 60, and 120.
The graph in Theorem~\ref{mainthm} evaluates to~32.

The point $L_1^3+\cdots+L_8^3$ at which we evaluated is defined by the following linear forms:
$L_1=(1,0,0,0,0)$,
$L_2=(0,1,0,0,0)$,
$L_3=(0,0,1,0,0)$,
$L_4=(0,0,0,1,0)$,
$L_5=(0,0,0,0,1)$,
$L_6=(1,1,1,1,1)$,
$L_7=(1,2,3,2,1)$,
$L_8=(1,2,1,1,2)$.
The gcd we divided by is $57600$, which is the value to which for example the 6th graph in our list evaluates.

\begin{Remark}
Using our calculation we can reprove Ottaviani's \cite{Ottaviani} result.
Indeed, it only remains to prove that the polynomial defined in Theorem~\ref{mainthm} is irreducible.
We use plethysm calculations in order to prove this.
Invariants of $SL_5$ correspond to rectangular Young diagrams (5 rows of equal lengths).
Since Young diagrams of irreducible representations in degree $d$ have exactly $d$ boxes,
invariants of $SL_5$ can only appear in degrees $d$ that are multiples of 5.
We used the \textsc{Schur} software, available online at \verb+http://sourceforge.net/projects/schur+,
to determine that there is no $SL_5$ invariant in degree 5
(which is very easy to see using the symbolic method),
but there is exactly one (up to scale) in degree 10 and exactly one (up to scale) in degree 15.
Assume for the sake of contradiction that the degree 15 invariant is a reducible polynomial.
It follows from the connectedness of $SL_5$ that all its factors are
invariant (the same fact was used in~\cite[p. 105]{Ottaviani}).
We thus get a nonzero degree 15 polynomial with  at least two factors, each of which has at least degree 10: a contradiction.
\end{Remark}

\begin{Remark}
It is quite remarkable that there is only one invariant in degree 15. As a result, {\em any} graph would give this
invariant provided the associated polynomial is not identically zero.
\end{Remark}

\section{A human proof of the chromatic property for the most symmetric graph}

\begin{figure}
\begin{center}
\begin{tikzpicture}[scale=0.80]
\tikzstyle{vertex}=[circle,draw=black, fill=white, inner sep = 0.75mm]
\node [shape=regular polygon, regular polygon sides = 5, inner sep = 0.75mm, draw=black, label=45:{\tiny $7$}] (v7) at (-1,-1){};
\node [shape=regular polygon, regular polygon sides = 3, draw=black, inner sep=0.55mm,  label=135:{\tiny $12$}] (v12) at (1,-1) {};
\node [shape=regular polygon, regular polygon sides = 3, draw=black, inner sep=0.55mm,  label=225:{\tiny $11$}] (v11) at (1,1) {};
\node [shape=regular polygon, regular polygon sides = 5,  draw=black, inner sep = 0.75mm,  label=315:{\tiny $2$}] (v2) at (-1,1) {};
\node [shape=star, draw=black, inner sep=0.75mm,  label=270:{\tiny $5$}] (v5) at (-2,-2) {};
\node [shape=star, draw=black, inner sep=0.75mm, label=90:{\tiny $1$}] (v1) at (-2,2) {};
\node [shape=diamond, draw=black, inner sep=0.75mm,   label=90:{\tiny $3$}] (v3) at (2,2) {};
\node [shape=diamond, draw=black, inner sep=0.75mm,  label=270:{\tiny $8$}] (v8) at (2,-2) {};

\node [vertex,  label=90:{\tiny $0$}] (v0) at (-5,1) {};
\node [vertex,  label=180:{\tiny $4$}] (v4) at (-6,0) {};
\node [vertex,  label=270:{\tiny $6$}] (v6) at (-5,-1) {};

\node [shape=regular polygon, regular polygon sides = 3, draw=black, inner sep=0.55mm,  label=270:{\tiny $14$}] (v14) at (6,-1) {};
\node [shape=diamond, draw=black, inner sep=0.75mm, label=90:{\tiny $9$}] (v9) at (6,1) {};
\node [shape=regular polygon, regular polygon sides=5, draw=black, inner sep=0.75mm,  label=90:{\tiny $10$}] (v10) at (8,1) {};
\node [shape=star, inner sep = 0.75mm,  draw=black, label=270:{\tiny $13$}] (v13) at (8,-1) {};

\draw (v0)--(v4)--(v6)--(v0);

\draw (v9)--(v10)--(v13)--(v14);
\draw (v13)--(v9)--(v14)--(v10);

\draw (v1)--(v3)--(v8)--(v5)--(v1);
\draw (v2)--(v11)--(v12)--(v7)--(v2);

\draw (v1)--(v2)--(v3)--(v11)--(v1);

\draw (v5)--(v7)--(v8)--(v12)--(v5);

\node (cj) at (-3.5,0) {\footnotesize all edges};
\node (dj) at (4,0) {\footnotesize no same symbol edges};

\end{tikzpicture}
\end{center}
\caption{The collinearity graph}
\label{fig:coll}
\end{figure}
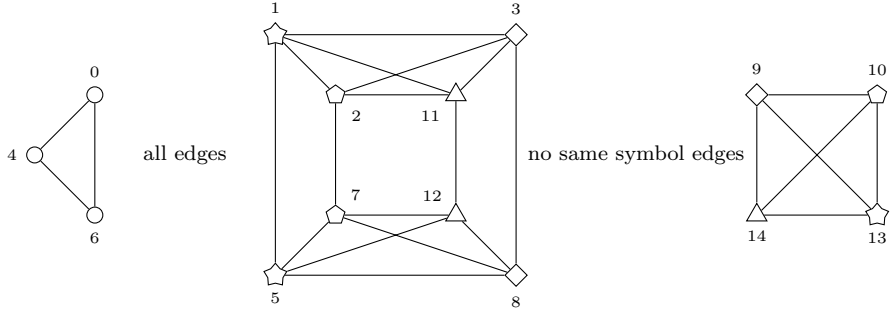

We consider the collinearity graph of the hypergraph from Theorem \ref{mainthm} and show ``by hand'' that it is not 7-colorable.
The automorphism group of the collinearity graph has three orbits, which induce the three subgraphs shown in Figure~\ref{fig:coll}, namely 
a triangle, a quartic graph on 8 vertices that can easily be identified as the complement of the 3-cube $Q_3$, and a 4-clique.
The edges between the three orbits are not  shown in the diagram, but between
the triangle and the complement of the cube all edges are present, while between the complement of the cube and the 4-clique, all edges
are present {\em except those joining matching symbols}. For example, vertex $13$, drawn with a star, is adjacent to all the vertices in the
complement of the cube except vertices $1$ and $5$ each of which is drawn by a star. There are no further edges. This graph has 
automorphism group $C_2 \times S_3 \times S_4$ where the $C_2$ exchanges the like-symbols in the complement of the cube, the $S_3$ arbitrarily
permutes the vertices in the triangle and the $S_4$ arbitrarily permutes the vertices in the $4$-clique (necessarily inducing
the same permutation on the symbols in each half of the complement of the cube). Notice that the automorphism group of the collinearity 
graph is considerably larger (size 288)
than the automorphism group of the hypergraph (size 12) because some permutations of the points preserve 
collinearity even though they do not preserve hyperedges. 

\bigskip

\newcommand{\cliq}[5]{
\draw (v#1)--(v#2);
\draw (v#1)--(v#3);
\draw (v#1)--(v#4);
\draw (v#1)--(v#5);
\draw (v#2)--(v#3);
\draw (v#2)--(v#4);
\draw (v#2)--(v#5);
\draw (v#3)--(v#4);
\draw (v#3)--(v#5);
\draw (v#4)--(v#5);
}

Now, suppose for a contradiction that there is a $7$-coloring of this graph. First notice that the triangle $\{0,4,6\}$ together with
the ``top-half'' of the complement of the cube, $\{1,2,3,11\}$ forms a $7$-clique in this graph, and hence it  uses all $7$ colors.
Similarly, $\{0,4,6\}$ together
with the ``bottom-half'' of the complement of the cube forms a $7$-clique which also uses all $7$ colors. Thus the sets of vertices
$\{1,2,3,11\}$ and $\{5,7,8,12\}$ use the same set of four colors. Now consider a vertex in the $4$-clique, which is adjacent to six vertices
in the complement of the cube, three from each half. As each pair of vertices with the same symbol is colored differently
in the top and bottom halves, this set of six vertices collectively use all four of the colors. Therefore the only colors available to
color the $4$-clique are the three colors that were assigned to the triangle, which is the required contradiction. 

\section{Other situations}

Let us call a triple $(k,d,n)$ {\em chromatic} if it satisfies all of the following properties:
\begin{enumerate}
\item
The secant variety $\sigma_k(v_d(\mathbb{P}^n))$ is of codimension 1.
\item
The defining equation of this secant variety can be expressed by a single graph (with the appropriate vertex degrees
following the same idea as in \S\ref{introsec}).
\item
There exists such a graph with chromatic number $k+1$.
\end{enumerate}

Note that by $SL_{n+1}$-invariance of the defining equation, provided Condition (1) holds, and by the First Fundamental Theorem
of Invariant Theory for $SL_{n+1}$, one can a priori express the defining polynomial as a linear combination
of graphs as in \S\ref{introsec}. However, there is no reason for only one graph to do the job as required in Condition (2).
If such a graph exists, there is also no reason for there being a simple coloring explanation for its vanishing on the secant variety.
There could, for instance, be some subtle cancellations between the contributions of proper $k$-colorings.
Nevertheless, our main result is that $(7,3,4)$ is chromatic. One may wonder if this is a general phenomenon when the secant variety
has codimension 1.
Triples satisfying the codimension 1 condition fall in two classes, according to the Alexander-Hirschowitz Theorem:
\begin{enumerate}
\item
AH-exceptional: $(n,2,n)$, for $n\ge 2$, as well as $(5,4,2)$, $(9,4,3)$, $(14,4,4)$, and $(7,3,4)$.
\item
AH-ordinary: $(k,d,n)$ which satisfy the Diophantine equation $\binom{n+d}{d}=(n+1)k+1$.
\end{enumerate}

The determinant of a quadratic form in $n+1$ variables
has symbolic expression
\[
\frac{1}{(n+1)!} (a^{(1)}a^{(2)}\cdots a^{(n+1)})^2
\]
where the superscripts $(i)$ label different symbolic letters (which themselves stand for vectors of $n+1$ variables).
The corresponding hypergraph on a set of vertices $\mathcal{A}$ with $n+1$ elements only has two hyperedges
or $(n+1)$-blocks both equal to $\mathcal{A}$
(here we have a repeated hyperedge contrary to the result of the search in \S\ref{graphsec}). Clearly such a triple
$(n,2,n)$ is chromatic.

For binary forms, one has that the AH-ordinary triples $(k,2k,1)$ are chromatic.
Indeed, the defining polynomial or classical catalecticant has symbolic expression
\[
\prod_{1\le i<j\le k+1}(a^{(i)}a^{(j)})^2\ .
\]
The corresponding hypergraph, now an ordinary (multi-)graph, is the complete graph
on $k+1$ elements with doubled edges.

Remark \ref{Sinvremark} also shows that the AH-ordinary triple $(3,3,2)$ is chromatic.
As for the AH-exceptional cases besides $(7,3,4)$ which has just been settled, one has the following results.
The invariant with symbolic expression
\[
(abc)^2(ade)(adf)(bdf)(bef)(cde)(cef)
\]
provides the equation for the $(5,4,2)$ case. An equation for this secant variety of degree 6
was known since Clebsch (see~\cite[p. 1234]{BrambillaO}).
The above symbolic formula may however be new. It is easy to see ``by hand'' that the collinearity graph
obtained from the hypergraph encoded in the above expression is the complete graph on 6 vertices.
The nonvanishing has been checked by computer
using the same method as in \S\ref{nonvansec}. Thus $(5,4,2)$ is chromatic.
Note that the symbolic expression given in~\cite[\S4]{Chipalkatti} (and therein denoted by $\Phi$(6,0,0))
does not produce a collinearity graph which is a $6$-clique.

Similarly, the invariant with symbolic formula
\[
(abcd)(abgj)(aefg)(afhi)(bdef)(behi)(cdgh)(ceij)(cfhj)(dgij)
\]
provides an equation in the $(9,4,3)$ case.
The corresponding collinearity graph is a $10$-clique.
The nonvanishing was checked by computer as in \S\ref{nonvansec}
and therefore $(9,4,3)$ is chromatic.

Unfortunately, $(14,4,4)$ is {\em not} chromatic.
In this case the degree of the equation is 15 (see~\cite[p. 1234]{BrambillaO})
and one would thus look for a $5$-uniform $4$-regular hypergraph on a set $\mathcal{A}$ with $15$ elements (and thus with
twelve $5$-blocks)
whose collinearity graph is a $15$-clique.
Such a hypergraph cannot exist.  Let $C(m,5,2)$ be defined as the smallest number of $5$-subsets (or blocks) of an $m$-set (of points)
such that every $2$-set of points is ``covered'' by at least one block (this is also called a covering number in the design theory literature).
Interestingly, for $m\equiv 3\pmod{4}$ the covering numbers $C(m,5,2)$ have all been determined~\cite{MillsM}
(see also~\cite[Ch. VI.11]{ColbournD}).
Let 
\[
B(m)=\left\lceil
\frac{m\left\lceil \frac{m-1}{4}\right\rceil}{5}
\right\rceil
\]
which in particular gives $B(15)=12$.
The main result in~\cite{MillsM} is that for all such $m$'s one has $C(m,5,2)=B(m)$ {\em except}
for $m=15$ where $C(15)=B(15)+1=13$ (this fact has been observed in~\cite{StantonKM}).
Quite remarkably, our case of interest turns out to be an exceptional
case from the point of view of design theory. Thus, even relaxing the regularity conditions, we are one $5$-block short
from being able to construct the needed $15$-clique.

\bigskip
\noindent{\bf Acknowledgements:}
{\small
A. A. thanks Jaydeep Chipalkatti and Achilleas Pitsillides for their help at an early stage of this project. He also
thanks Andrew Obus for a useful discussion.
C. I. acknowledges the Texas A\&M University Brazos HPC cluster that contributed to the research
reported here. $\langle$brazos.tamu.edu$\rangle$.
G. R.'s research was partially supported by Australian Research Council Discovery Project grants.
}

\end{document}